\documentclass[a4paper,12pt]{amsart}
\usepackage{amssymb, amsmath}
\usepackage[curve]{xypic}

\providecommand{\abs}[1]{\lvert#1\rvert}

\providecommand{\Ob}{\textnormal{Ob}}

\providecommand{\cpt}{\textnormal{cpt}}

\begin{document}

\title{Geometric homology revisited}
\author{Fabio Ferrari Ruffino}
\address{ICMC - Universidade de S\~ao Paulo, Avenida Trabalhador s\~ao-carlense 400, 13566-590 - S\~ao Carlos - SP, Brasil}
\email{ferrariruffino@gmail.com}
\thanks{The author was supported by FAPESP (Funda\c{c}\~ao de Amparo \`a Pesquisa do Estado de S\~ao Paulo).}

\begin{abstract}
Given a cohomology theory $h^{\bullet}$, there is a well-known abstract way to define the dual homology theory $h_{\bullet}$, using the theory of spectra. In \cite{Jakob} the author provides a more geometric construction of $h_{\bullet}$, using a generalization of the bordism groups. Such a generalization involves in its definition the vector bundle modification, which is a particular case of the Gysin map. In this paper we provide a more natural variant of that construction, which replaces the vector bundle modification with the Gysin map itself, which is the natural push-forward in cohomology. We prove that the two constructions are equivalent.
\end{abstract}

\maketitle

\newtheorem{Theorem}{Theorem}[section]
\newtheorem{Lemma}[Theorem]{Lemma}
\newtheorem{Corollary}[Theorem]{Corollary}
\newtheorem{Rmk}[Theorem]{Remark}
\newtheorem{Def}{Definition}[section]
\newtheorem{ThmDef}[Theorem]{Theorem - Definition}

\section{Introduction}

Given a cohomology theory $h^{\bullet}$, there is a well-known abstract way to define the dual homology theory $h_{\bullet}$, using the theory of spectra. In particular, if $h^{\bullet}$ is representable via a spectrum $E = \{E_{n}, e_{n}, \varepsilon_{n}\}_{n \in \mathbb{Z}}$, for $e_{n}$ the marked point of $E_{n}$ and $\varepsilon_{n}: \Sigma E_{n} \rightarrow E_{n+1}$ the structure map, one can define on a space with marked point $(X, x_{0})$ \cite{Rudyak}:
	\[h_{n}(X, x_{0}) := \pi_{n}(E \wedge X).
\]
In \cite{Jakob} the author provides a more geometric construction of $h_{\bullet}$, using a generalization of the bordism groups. In particular, he shows that, for a given pair $(X,A)$, a generator of $h_{n}(X,A)$ can be represented by a triple $(M, \alpha, f)$, where $M$ is a compact $h^{\bullet}$-manifold with boundary of dimension $n+q$, $\alpha \in h^{q}(M)$ and $f: (M, \partial M) \rightarrow (X,A)$ is a continuous function. On such triples one must impose a suitable equivalence relation, which is defined via the natural generalization of the bordism equivalence relation and via the notion of vector bundle modification. In this paper we provide a variant of that construction, which seems to be more natural. In particular, we replace the notion of vector bundle modification with the Gysin map, the latter being the natural push-forward in cohomology. The vector bundle modification is just a particular case, which holds when the underlying map is a section of a sphere bundle. We prove that the two constructions are equivalent, since there is a natural isomorphism between the geometric homology groups defined in \cite{Jakob} and their variant defined in the present paper.

The paper is organized as follows. In section \ref{Preliminaries} we recall the definition of the Gysin map, even for manifolds with boundary, and the geometric construction of the homology groups provided in \cite{Jakob}. In section \ref{DualRevisited} we introduce the variant of the geometric construction we discussed above, and in section \ref{Equivalence} we prove that the two constructions are equivalent.

\section{Preliminaries}\label{Preliminaries}

We call $\mathcal{FCW}_{2}$ the category of pairs of spaces having the homotopy type of finite CW-complexes. Let $h^{\bullet}$ be a multiplicative cohomology theory on $\mathcal{FCW}_{2}$. We recall the construction of the Gysin map for smooth maps between differentiable manifolds with boundary (v.\ \cite{Karoubi, FR} for manifolds without boundary).

\subsection{Gysin map}

Let $h^{\bullet}$ be a cohomology theory on $\mathcal{FCW}_{2}$. A smooth $h^{\bullet}$-manifold is a smooth manifold with an $h^{\bullet}$-orientation, the latter being a Thom class of its tangent bundle or, equivalently, of its stable normal bundle. Given two compact smooth $h^{\bullet}$-manifolds with boundary $X$ and $Y$ and a map $f: (Y, \partial Y) \rightarrow (X, \partial X)$, we can define the Gysin map $f_{!}: h^{\bullet}(Y) \rightarrow h^{\bullet + \dim\,X - \dim\,Y}(X)$ as:
\begin{equation}\label{GysinLefschetz}
	f_{!}(\alpha) := L_{X}^{-1} f_{*} L_{Y} (\alpha)
\end{equation}
where:
\begin{equation}\label{Lefschetz}
	L_{X}: h^{\bullet}(X) \rightarrow h_{\dim X - \bullet}(X, \partial X)
\end{equation}
is the Lefschetz duality \cite{Switzer}. The problem of this definition is that it involves the homology groups, which we have to define, therefore we need a construction involving only the cohomology groups. When $f$ is a neat map, one can define the Gysin map in a way similar to the one shown in \cite{Karoubi}, pp.\ 230-234, about topological K-theory on manifolds without boundary. Then the definition can be easily extended to any map between $h^{\bullet}$-manifolds.

We start with embeddings. We call $\mathbb{R}^{n}_{+} := \{(x_{1}, \ldots, x_{n}) \in \mathbb{R}^{n} \,\vert\, x_{n} \geq 0\}$. Given a manifold $X$ and a point $x \in \partial X$, by definition there exists a chart $(U, \varphi)$ of $X$ in $x$, with $U \subset X$ open and $\varphi: U \rightarrow \mathbb{R}^{n}_{+}$, such that $\varphi(x) = 0$ and $\varphi(\partial X \cap U) = (\mathbb{R}^{n-1} \times \{0\}) \cap \varphi(U)$. We call such a chart \emph{boundary chart}. For $m \leq n-1$, we call $\overline{\mathbb{R}}^{m}$ the subspace of $\mathbb{R}^{n}$ containing those vectors whose first $n-m$ components are vanishing, i.e.\ $\overline{\mathbb{R}}^{m} = \{0\} \times \mathbb{R}^{m} \subset \mathbb{R}^{n}$, and $\overline{\mathbb{R}}^{m}_{+} := \overline{\mathbb{R}}^{m} \cap \mathbb{R}^{n}_{+}$.
\begin{Def} An embedding of manifolds $i: (Y, \partial Y) \hookrightarrow (X, \partial X)$, where $\dim Y = m$, is \emph{neat} \cite{Hirsh, Kosinski} if:
\begin{itemize}
	\item $i(\partial Y) = i(Y) \cap \partial X$;
	\item for every $y \in \partial Y$ there exists a boundary chart $(U, \varphi)$ of $X$ in $i(y)$ such that $U \cap i(Y) = \varphi^{-1}(\overline{\mathbb{R}}^{m}_{+})$.
\end{itemize}
\end{Def}
The importance of neat embeddings in this context relies on fact that the properties of tubular neighborhoods are similar to the ones holding for manifolds without boundary.
\begin{Def} Let $(Y, \partial Y)$ be a neat submanifold of $(X, \partial X)$. A tubular neighborhood $U$ of $Y$ in $X$ is \emph{neat} \cite{Kosinski} if $U \cap \partial X$ is a tubular neighborhood of $\partial Y$ in $\partial X$.
\end{Def}
\begin{Theorem} If $(Y, \partial Y)$ is a neat submanifold of $(X, \partial X)$, there exists a neat tubular neighborhood of $Y$ in $X$ and it is unique up to isotopy.
\end{Theorem}
The proof can be found in \cite[Chapter~4.6]{Hirsh} and in \cite[Chapter~III.4]{Kosinski}. Let $i: (Y, \partial Y) \hookrightarrow (X, \partial X)$ be a neat embedding of smooth compact manifolds of codimension $n$, such that the normal bundle $N_{Y}X$ is $h^{\bullet}$-orientable. Let $U$ be a tubular neighborhood of $Y$ in $X$, and $\varphi_{U}: U \rightarrow N_{Y}X$ a homeomorphism, which exists by definition. The map:
	\[i_{!}: h^{\bullet}(Y) \rightarrow h^{\bullet + n}(X)
\]
is defined in the following way:
\begin{itemize}
	\item we apply the Thom isomorphism $T: h^{\bullet}(Y) \rightarrow h^{\bullet + n}_{\cpt}(N_{Y}X) = \tilde{h}^{\bullet + n}(N_{Y}X^{+})$, for $N_{Y}X^{+}$ the one-point compactification of $N_{Y}X$;
	\item we extend $\varphi_{U}$ to $\varphi_{U}^{+}: U^{+} \rightarrow N_{Y}X^{+}$ in the natural way and apply $(\varphi_{U}^{+})^{*}: h^{\bullet}_{\cpt}(N_{Y}X) \rightarrow h^{\bullet}_{\cpt}(U)$;
	\item considering the natural map $\psi: X \rightarrow U^{+}$ given by:
	\[\psi(x) = \left\{\begin{array}{ll}
	x & \text{if } x \in U \\
	\infty & \text{if } x \in X \setminus U
	\end{array}\right.
\]
we apply $\psi^{*}: \tilde{h}^{\bullet}(U^{+}) \rightarrow \tilde{h}^{\bullet}(X)$.
\end{itemize} 
Summarizing:
\begin{equation}\label{GysinMap}
	i_{!}(\alpha) := \psi^{*} \circ (\varphi_{U}^{+})^{*} \circ T(\alpha).
\end{equation}
We now define the Gysin map associated to a generic neat map $f: (Y, \partial Y) \rightarrow (X, \partial X)$, not necessarily an embedding.
\begin{Def} A smooth map $f: (Y, \partial Y) \rightarrow (X, \partial X)$ is \emph{neat} (v.\ \cite[Appendix~C]{HS} and references therein) if:
\begin{itemize}
	\item $f^{-1}(\partial X) = \partial Y$;
	\item for every $y \in \partial Y$, the map $df_{y}: T_{y}Y/T_{y}\partial Y \rightarrow T_{f(y)}X/T_{f(y)}\partial X$ is an isomorphism.
\end{itemize}
\end{Def}
If $f$ is an embedding this definition is equivalent to the previous one. In the case of manifolds without boundary, in order to construct the Gysin map one considers an embedding $j: Y \rightarrow \mathbb{R}^{N}$, and the embedding $(f, j): Y \rightarrow X \times \mathbb{R}^{N}$. This does not apply to manifolds with boundary, since $j$ is not a neat map, and, if we consider $\mathbb{R}^{N}_{+}$ instead of $\mathbb{R}^{N}$, then it is more complicated to define the integration map. Anyway, a similar construction is possible thanks to the following theorem (v.\ \cite[Appendix~C]{HS} and references therein).
\begin{Theorem} Let $f: (Y, \partial Y) \rightarrow (X, \partial X)$ be a neat map. Then there exists a neat embedding $\iota: (Y, \partial Y) \rightarrow (X \times \mathbb{R}^{N}, \partial X \times \mathbb{R}^{N})$, stably unique up to isotopy, such that $f = \pi_{X} \circ \iota$ for $\pi_{X}: X \times \mathbb{R}^{N} \rightarrow X$ the projection.
\end{Theorem}
Therefore we consider the Gysin map:
	\[\iota_{!}: h^{\bullet}(Y) \rightarrow h_{\cpt}^{\bullet + (N + \dim X - \dim Y)}(X \times \mathbb{R}^{N})
\]
as previously defined, followed by the integration map:
\begin{equation}\label{IntegrationMap}
	\int_{\mathbb{R}^{N}}: \; h^{\bullet+N}_{\cpt}(X \times \mathbb{R}^{N}) \rightarrow h^{\bullet}(X)
\end{equation}
defined in the following way:
\begin{itemize}
	\item $h^{\bullet+N}_{\cpt}(X \times \mathbb{R}^{N}) = \tilde{h}^{\bullet+N}((X \times \mathbb{R}^{N})^{+}) \simeq \tilde{h}^{\bullet+N}(\Sigma^{N}(X_{+}))$, for $X_{+} = X \sqcup \{\infty\}$;
	\item we apply the suspension isomorphism $\tilde{h}^{\bullet+N}(\Sigma^{N}(X_{+})) \simeq \tilde{h}^{\bullet}(X_{+}) \simeq h^{\bullet}(X)$.
\end{itemize}
Summarizing:
	\[f_{!}(\alpha) := \int_{\mathbb{R}^{N}} \iota_{!}(\alpha).
\]
In order to prove that the Gysin map so defined does not depend on the choices involved in the construction (the tubular neighborhood $U$, the diffeomorphism $\varphi_{U}$ with the normal bundle, the embedding $\iota$), the proof in \cite{Karoubi} applies also to the case of manifolds with boundary. In fact, the independence from the tubular neighborhood and the associated diffeomorphism is a consequence of the unicity up to isotopy (in particular, homotopy) of such a neighborhood. Also for what concerns the embedding $\iota$, the proof of \cite{Karoubi}, Prop.\ 5.24 p.\ 233, applies. In particular, $f_{!}$ only depends on the homotopy class of $f$.

If $f: (Y, \partial Y) \rightarrow (X, \partial X)$ is a generic map, not necessarily neat, we can define $f_{!}$ via the following lemma.
\begin{Lemma}\label{HomotopicNeat} Any smooth map $f: (Y, \partial Y) \rightarrow (X, \partial X)$ between compact manifolds is homotopic to a neat map relatively to $\partial Y$.
\end{Lemma}
\paragraph{Proof:} We can choose two collar neighborhoods $U$ of $\partial Y$ and $V$ of $\partial X$, such that $f(U) \subset V$. Hence we think of $f\vert_{U}$ as a map from $\partial Y \times [0, 1)$ to $\partial X \times [0, 1)$. We consider the following homotopy: $F_{t}(y, u) = (\pi_{\partial X}f(y, u), (1-t)\pi_{[0,1)}f(y,u) + tu)$. It follows that $F_{0} = f$ and $F_{1}(y, u) = (\pi_{\partial X}f(y, u), u)$, and the latter is neat. Gluing $F_{t}$ and $f\vert_{X \setminus U}$ via a bump function, we get a complete homotopy. $\square$ \\

Since the Gysin map $f_{!}$, for $f$ neat, only depends on the homotopy class of $f$, we can define it for a generic map $f$, simply considering any neat function homotopic to it. The Gysin map commutes with the restrictions to the boundaries up to a sign, as the following theorem shows.
\begin{Theorem}\label{GysinCommBoundary} Let $f: (Y, \partial Y) \rightarrow (X, \partial X)$ be a map between $h^{\bullet}$-oriented smooth manifolds, and $f': \partial Y \rightarrow \partial X$ the restriction to the boundaries. Then, for $\alpha \in h^{\bullet}(X)$:
	\[f'_{!}(\alpha\vert_{\partial Y}) = (-1)^{\dim X - \dim Y}(f_{!}\alpha)\vert_{\partial X}
\]
where the orientations on the boundaries are naturally induced from the ones of $X$ and $Y$.
\end{Theorem}
\paragraph{Proof:} It is enough to prove the statement for embeddings, since the integration over $\mathbb{R}^{N}$, which is actually the suspension isomorphism, commutes with the restrictions. Therefore, let us suppose that $f$ is a neat embedding, and that $N_{Y}X$ is a neat tubular neighborhood. Then $N_{\partial Y}\partial X = N_{Y}X\vert_{\partial Y}$, but the orientation induced by the ones of $\partial Y$ and $\partial X$ on $N_{\partial Y}\partial X$ differs from the restriction of the one induced by $Y$ and $X$ by a factor $(-1)^{\dim X - \dim Y}$. Therefore, if:
	\[T: h^{\bullet}(Y) \rightarrow h^{\bullet + \dim X - \dim Y}_{\cpt}(N_{Y}X), \qquad T': h^{\bullet}(\partial Y) \rightarrow h^{\bullet + \dim X - \dim Y}_{\cpt}(N_{\partial Y}\partial X)
\]
are the Thom isomorphisms, it follows that:
	\[T'(\alpha\vert_{\partial Y}) = (-1)^{\dim X - \dim Y}T(\alpha)\vert_{N_{\partial Y}\partial X}.
\]
Then, since the pull-backs commute with the restrictions, the thesis follows. If $f$ is a generic embedding, not necessarily neat, since $f$ is homotopic to a neat embedding relatively to the boundary, then $f'$ remains unchanged under the homotopy and the thesis follows. $\square$

\subsection{Geometric homology}

We recall the geometric definition of the homology theory dual to a given cohomology theory $h^{\bullet}$, as defined in \cite{Jakob}. Let $M$ be a paracompact space, $\pi_{V}: V \rightarrow M$ a real vector bundle of rank $r+1$ with metric and $h^{\bullet}$-orientation, and $\sigma: M \rightarrow V$ a section of norm $1$. Then, $\sigma$ induces an isomorphism $V \simeq E \oplus 1$, for $\pi_{E}: E \rightarrow M$ a vector bundle of rank $r$ with metric, such that $\sigma(m) \simeq (0, 1)_{m}$. We identify $V$ with $E \oplus 1$. The unit sphere bundle $SV$ of $V$ can be thought of as the union of two disc bundles, the two hemispheres, joined on the equator: the two disc bundles are isomorphic to the unit disc bundle of $E$, therefore we call them $D^{+}E$ and $D^{-}E$, while the bundle of the equators is isomorphic to the sphere bundle of $E$, which we call $SE$. Moreover, the north pole $(0, 1)_{m}$ of $D^{+}E_{m}$ is $\sigma(m)$, therefore $D^{+}E$ is a tubular neighborhood of the image of $\sigma$. There is a natural map:
\begin{equation}\label{VBM}
	\sigma_{!}: h^{\bullet}(M) \rightarrow h^{\bullet+r}(SV)
\end{equation}
defined in the following way:
\begin{itemize}
	\item we apply the Thom isomorphism $T: h^{\bullet}(M) \rightarrow h^{\bullet+r}(E^{+}) \simeq h^{\bullet+r}(D^{+}E, SE)$;
	\item by excision $h^{\bullet+r}(D^{+}E, SE) \simeq h^{\bullet+r}(SV, D^{-}E)$;
	\item from the inclusion of couples $(SV, \emptyset) \subset (SV, D^{-}E)$ we get a map $h^{\bullet+r}(SV, D^{-}E) \rightarrow h^{\bullet+r}(SV)$.
\end{itemize}
The map \eqref{VBM} coincides with the Gysin map associated to $\sigma$ \cite{Jakob}.
\begin{Def}\label{CyclesJakob} For $(X, A) \in \Ob(\mathcal{FCW}_{2})$ and $n \in \mathbb{Z}$ fixed, we consider the quadruples $(M, u, \alpha, f)$ where:
\begin{itemize}
	\item $(M, u)$ a smooth compact $h^{\bullet}$-manifold, possibly with boundary, whose connected components $\{M_{i}\}$ have dimension $n+q_{i}$, with $q_{i}$ arbitrary; we think of $u$ as a Thom class of the tangent bundle;
	\item $\alpha \in h^{\bullet}(M)$, such that $\alpha\vert_{M_{i}} \in h^{q_{i}}(M)$;
	\item $f: (M, \partial M) \rightarrow (X, A)$ is a map.
\end{itemize}
Two quadruples $(M, u, \alpha, f)$ and $(N, v, \beta, g)$ are equivalent if there exists an orientation-preserving diffeomorphism $F: (M, u) \rightarrow (N, v)$ such that $f = g \circ F$ and $\alpha = F^{*}\beta$. The \emph{group of $n$-cycles} $C_{n}(X, A)$ is the free abelian group generated by equivalence classes of such quadruples.
\end{Def}
We now consider the group $G_{n}(X, A)$ defined as the quotient of $C_{n}(X, A)$ by the subgroup generated by elements of the form:
\begin{itemize}
	\item $[(M, u, \alpha, f)] - [(M_{1}, u\vert_{M_{1}}, \alpha\vert_{M_{1}}, f\vert_{M_{1}})] - [(M_{2}, u\vert_{M_{2}}, \alpha\vert_{M_{2}}, f\vert_{M_{2}})]$, for $M = M_{1} \sqcup M_{2}$;
	\item $[(M, u, \alpha + \beta, f)] - [(M, u, \alpha, f)] - [(M, u, \beta, f)]$.
\end{itemize}
Moreover, we define the subgroup $U_{n}(X, A) \leq G_{n}(X, A)$ as the one generated by elements:
\begin{itemize}
	\item $[(M, u, \alpha, f)] - [(S(E \oplus 1), \tilde{u}, \sigma_{!}\alpha, f \circ \pi)]$, where $S(E \oplus 1)$ is the sphere bundle induced by an $h^{\bullet}$-oriented vector bundle\footnote{The vector bundle $E$ may have different rank on different connected components of $M$.} $E \rightarrow M$ with metric, $\tilde{u}$ is the orientation canonically induced on $S(E \oplus 1)$ as a manifold from $u$ and the orientation of $E$, $\sigma: M \rightarrow E \oplus 1$ is the section $\sigma(m) = (0, 1)_{m}$ and $\sigma_{!}$ is the vector bundle modification \eqref{VBM};
	\item $[(M, u, \alpha, f)]$ such that there exists $[(W, U, A, F)] \in G_{n+1}(X, X)$ such that $M \subset \partial W$ is a regularly embedded submanifold of codimension $0$ and $F(\partial W \setminus M) \subset A$, $U = u\vert_{M}$, $\alpha = A\vert_{M}$, $f = F\vert_{M}$.
\end{itemize}
Finally:
\begin{Def} The geometric homology groups are defined as $h_{n}(X, A) := G_{n}(X, A) / U_{n}(X, A)$.
\end{Def}

\section{Geometric homology revisited}\label{DualRevisited}

We now redefine the homology groups using only the Gysin map instead of the vector bundle modification. We also define cycles and boundaries in a slightly different way.
\begin{Def} On a pair $(X, A) \in \mathcal{FCW}_{2}$, we define the group of \emph{$n$-precycles of $h_{\bullet}$} as the free abelian group generated by the quadruples $(M, u, \alpha, f)$, for:
\begin{itemize}
	\item $(M, u)$ a smooth compact $h^{\bullet}$-manifold, possibly with boundary, whose connected components $\{M_{i}\}$ have dimension $n+q_{i}$, with $q_{i}$ arbitrary;
	\item $\alpha \in h^{\bullet}(M)$, such that $\alpha\vert_{M_{i}} \in h^{q_{i}}(M)$;
	\item $f: (M, \partial M) \rightarrow (X, A)$ a continuous map.
\end{itemize}
\end{Def}
Contrary to definition \ref{CyclesJakob}, we do not quotient out with respect to orientation-preserving diffeomorphisms, since it will turn out not to be necessary. We define cycles and boundaries in the following way.
\begin{Def} The group of \emph{$n$-cycles of $h_{\bullet}$}, denoted by $z_{n}(X, A)$, is the quotient of the group of $n$-precycles by the free subgroup generated by elements of the form:
\begin{itemize}
	\item $(M, u, \alpha, f) - (M_{1}, u\vert_{M_{1}}, \alpha\vert_{M_{1}}, f\vert_{M_{1}}) - (M_{2}, u\vert_{M_{2}}, \alpha\vert_{M_{2}}, f\vert_{M_{2}})$, for $M = M_{1} \sqcup M_{2}$;
	\item $(M, u, \alpha + \beta, f) - (M, u, \alpha, f) - (M, u, \beta, f)$;
	\item $(M, u, \varphi_{!}\alpha, f) - (N, v, \alpha, f \circ \varphi)$ for $\varphi: (N, \partial N) \rightarrow (M, \partial M)$ a map.
\end{itemize}
\end{Def}
The use of the Gysin map in this definition is more natural than the one of the vector bundle modification, which is just a particular case, while the Gysin map is the natural push-forward, defined for any $\varphi: (N, \partial N) \rightarrow (M, \partial M)$. Moreover, it is not necessary to explicitly deal with diffeomorphisms: in fact, if $\varphi: (M, u) \rightarrow (N, v)$ is an orientation preserving diffeomorphism, it is trivial to show from the definition that $\varphi_{!} = (\varphi^{-1})^{*}$, therefore the quotient in definition \ref{CyclesJakob} is just another particular case of the Gysin map.
\begin{Def} The group of \emph{$n$-boundaries of $h_{\bullet}$}, denoted by $b_{n}(X, A)$, is the subgroup of $z_{n}(X, A)$ containing the cycles which are representable by a precycle $(M, u, \alpha, f)$ such that there exits a precycle $(W, U, A, F)$ of $(X, X)$ such that:
\begin{itemize}
	\item $M \subset \partial W$ is a regularly embedded submanifold of codimension $0$;
	\item $F(\partial W \setminus M) \subset A$;
	\item $U = u\vert_{M}$, $\alpha = A\vert_{M}$, $f = F\vert_{M}$.
\end{itemize}
\end{Def}
Of course we define:
	\[h_{n}(X,A) := z_{n}(X,A) / b_{n}(X,A).
\]
For $g: (X, A) \rightarrow (Y, B)$ a map, the push-forward $g_{*}: h_{\bullet}(X,A) \rightarrow h_{\bullet}(Y, B)$ is naturally defined as $g_{*}[(M, u, \alpha, f)] = [(M, u, \alpha, f \circ g)]$, while the connecting homomorphism $\partial_{n}: h_{n}(X,A) \rightarrow h_{n-1}(A)$ is defined as:
	\[\partial_{n}[(M, u, \alpha, f)] = (-1)^{\abs{\alpha}}[(\partial M, u\vert_{\partial M}, \alpha\vert_{\partial M}, f\vert_{\partial M})]
\]
where $(-1)^{\abs{\alpha}}$ depends on the connected component and $u\vert_{\partial M}$ is the orientation naturally induced by $u$ on the boundary. It is well-defined thanks to theorem \ref{GysinCommBoundary}. The exterior product and the cap product are defined as in \cite{Jakob}.

\section{Equivalence}\label{Equivalence}

We call $h''_{n}(X, A)$ the geometric homology groups defined in \cite{Jakob}, $h'_{n}(X,A)$ the ones defined in the present paper and $h_{n}(X,A)$ the ones defined via spectra. There is a natural map:
\begin{equation}\label{EquivalenceMap}
\begin{split}
	\Psi_{n}(X, A):\; &h''_{n}(X, A) \rightarrow h'_{n}(X,A) \\
	& [(M, u, \alpha, f)] \rightarrow [(M, u, \alpha, f)],
\end{split}
\end{equation}
where of course the square brackets denote two different equivalence relations in the domain and in the codomain. It is easy to show that $\Psi_{n}$ is well-defined, since the quotient by diffeomorphisms and vector bundle modifications in the domain corresponds to quotient by the Gysin map in the codomain, therefore equivalent quadruples are sent to equivalent quadruples. Moreover, the quotient by boundaries, disjoint union of base manifolds and addition of cohomology classes are defined in the same way in the two cases.
\begin{Theorem}\label{EquivalenceTheorem} The maps $\Psi_{\bullet}$ defined in \eqref{EquivalenceMap} induce an equivalence of homology theories on $\mathcal{FCW}_{2}$ between $h''_{\bullet}$ and $h'_{\bullet}$.
\end{Theorem}
\paragraph{Proof:} $\Psi_{n}$ is a group homomorphism by construction, since it is defined on the generators of the free abelian group $C_{n}(X, A)$, and $h'_{n}(X, A)$ is defined quotienting out repeatedly $C_{n}(X, A)$. It is clearly surjective, since the elements of the form $[(M, u, \alpha, f)]$ are generators even for $h_{n}(X,A)$. Therefore, it remains to prove injectivity. Every element of $h''_{n}(X, A)$ can be written as a single class $[(M, u, \alpha, f)]$, since a sum of such classes can be reduced to a single one via the disjoint union of the base manifolds. Therefore, since $\Psi_{n}[(M, u, \varphi_{!}\alpha, f)] = \Psi_{n}[(N, v, \alpha, f \circ \varphi)]$ for $\varphi: (N, \partial N) \rightarrow (M, \partial M)$ a map, we must prove that $[(M, u, \varphi_{!}\alpha, f)] = [(N, v, \alpha, f \circ \varphi)]$ even in $h''_{n}(X,A)$. The equivalence between $h''_{n}(X,A)$ and $h_{n}(X,A)$ is given in \cite{Jakob} by the maps $\Phi_{\bullet}: h''_{\bullet}(X,A) \rightarrow h_{\bullet}(X,A)$ defined by $\Phi_{n}[(M, u, \alpha, f)] := f_{*}(\alpha \cap [M])$, for $[M]$ the fundamental class of $M$. By definition $\alpha \cap [M] = L_{M}(\alpha)$ for $L_{M}$ the Lefschetz duality \eqref{Lefschetz}. Hence, because of \eqref{GysinLefschetz}:
	\[\begin{split}
	& \Phi_{n}[(M, u, \varphi_{!}\alpha, f)] = f_{*} L_{M}(\varphi_{!}\alpha) \\
	& \Phi_{n}[(N, v, \alpha, f \circ \varphi)] = (f \circ \varphi)_{*} L_{N}(\alpha) = f_{*} \varphi_{*} L_{N}(\alpha) = f_{*} L_{M}(\varphi_{!}\alpha).
\end{split}\]
Since the map $\Phi_{n}$ is injective, it follow that $[(M, u, \varphi_{!}\alpha, f)] = [(N, v, \alpha, f \circ \varphi)]$ even in $h''_{n}(X,A)$. The fact that $\Psi_{\bullet}$ is a morphism of homology theories then follows from the fact that the boundary morphism and the push-forward are defined in the same way for $h'_{n}(X,A)$ and $h''_{n}(X,A)$. $\square$

\section*{Acknowledgements}

The author is financially supported by FAPESP (Funda\c{c}\~ao de Amparo \`a Pesquisa do Estado de S\~ao Paulo). We would like to thank Ryan Budney for having suggested the proof of lemma \ref{HomotopicNeat} on http://mathoverflow.net/.


\end{document}